\documentclass[12pt,reqno]{amsart}

\usepackage{epsfig}
\usepackage{amssymb,amsmath,amsthm,url}
\usepackage[a4paper,left=1in,top=1in,right=1in,bottom=1.1in]{geometry}

\usepackage{natbib}
\usepackage{graphicx}
\usepackage{amsfonts}
\usepackage{enumerate}
\usepackage[shortlabels]{enumitem}

\usepackage{IEEEtrantools}
\usepackage{algorithmic}
\usepackage{algorithm}
\usepackage{float}
\usepackage[colorinlistoftodos,prependcaption,textsize=tiny]{todonotes}
\usepackage{epsfig}
\usepackage{amssymb,amsmath,amsthm,url}
\usepackage{multirow}
\usepackage{booktabs}
\usepackage{natbib}
\usepackage{setspace}
\usepackage{dsfont}
\usepackage{amsmath}
\usepackage{amssymb}
\usepackage{amsthm}
\usepackage{upgreek}
\usepackage{mathrsfs}
\usepackage{makecell}
\usepackage{verbatim}
\usepackage[symbol]{footmisc}
\usepackage{lscape}
\usepackage{caption}
\usepackage{enumerate}
\usepackage{hyperref}
\usepackage{times,amsfonts,eucal}
\usepackage{times,amsfonts}
\usepackage{algorithm}
\usepackage{algorithmic}
\usepackage{graphicx,ifthen}
\usepackage{wrapfig}
\usepackage{latexsym}
\usepackage{color}
\usepackage{bm}
\usepackage{bbm}
\usepackage{srctex}

\newtheorem{theorem}{Theorem}[section]

\newtheorem{lemma}{Lemma}[section]
\newtheorem{assumption}{Assumption}
\newtheorem{remark}{Remark}

\begin{document}

\title[ ]{Consistency of binary segmentation for multiple change point estimation with functional data }

\author{Gregory Rice}
\address{Department of Statistics and Actuarial Science, University of Waterloo, Waterloo, ON, Canada}

\author{Chi Zhang}
\address{Department of Statistics and Actuarial Science, University of Waterloo, Waterloo, ON, Canada}

\begin{abstract}
For sequentially observed functional data exhibiting multiple change points in the mean function, we establish consistency results for the estimated number and locations of the change points based on the norm of the functional CUSUM process and standard binary segmentation. In addition to extending similar results in \cite{venkatraman:1992} and \cite{fryzl:2014:wildbin} for scalar data to the general Hilbert space setting, our main results are established without assuming the Gaussianity of the data, and under general linear process conditions on the model errors.
\end{abstract}
\maketitle
\section{Introduction}

Functional data analysis has emerged as a vibrant area of research in statistics over the past several decades, likely owing to the multitude of data collected, often at a high resolution, over a continuum. Such data can be viewed as discrete observations from functional data objects taking values in a function space. In a number of examples, functional data objects are obtained sequentially as a functional time series. We refer to \cite{HKbook} and \cite{ramsay:silverman:2002} for textbook length treatments of functional data analysis, and to \cite{bosq:2000} and \cite{hormann:kokoszka:2012} for reviews of functional time series analysis.

In the setting of functional time series, one often encounters series of curves that exhibit nonstationarity that appear as ``shocks" or structural changes in the data generating mechanism. A simple model for such data is a change point model, in which various features of the series are allowed to change at unknown points over the observation period. Recently a number of authors have studied change point models and methods to estimate them based on sequentially observed functional data, see \cite{aue:gabrys:horvath:kokoszka:2009}, \cite{berkes:gabrys:horvath:kokoszka:2009}, \cite{chiou:2019fdamulti}, \cite{aue:rice:sonmez:2018}, \cite{aston:kirch:2012AAS}, \cite{aston:kirch:2012JMVA}, \cite{bucchia:wendler:2017} and \cite{sharipov:wendler:2016}, who study changes point models for the mean structure, and \cite{Stoehr:2019}, \cite{Sharipov:2019}, and \cite{Dette2019} who consider changes in the higher order moment structure. Outside of \cite{chiou:2019fdamulti}, these papers generally consider the setting in which at most one change point is allowed in the model, and when this is not the case binary segmentation is proposed as a reasonable method to estimate and infer more than one change point.

Binary segmentation is an intuitive method to extend procedures for identifying a single change point to identify multiple changes, in which the sample is repeatedly segmented into two sub-samples based on estimates of a single change point until some stopping criterion is satisfied in order to identify further change points. \cite{vostrikova:1981} is usually credited with first having introduced the idea, and   \cite{venkatraman:1992} established the consistency of the procedure based on the standard CUSUM process for identifying and localizing changes in the mean of an independent and homoscedastic Gaussian sequence. Though a number of more modern techniques for identifying multiple changes points have been put forward, for example the PELT (\cite{killick:2012}) and MOSUM (\cite{eichinger:kirch:2018MOSUM} and \cite{hus:slaby:2001}) procedures, binary segmentation continues to be improved upon, see \cite{cho:fry:2012:multchange} and \cite{fryzl:2014:wildbin}, and exhibits strong and competitive performance, even in high-dimensional and non-stanard settings where more modern techniques have been less explored, see e.g. \cite{cho:fry:2015:highdimmult}.

The goal of this note is to establish the consistency of binary segmentation in detecting and localizing change points in the mean of a functional time series based on the norm of a functional analog of the CUSUM process. Beyond generalizing the results of \cite{venkatraman:1992} and the consistency results of \cite{fryzl:2014:wildbin} regarding binary segmentation to the general separable Hilbert space setting, 1) our main results do not require Gaussianity of the observations, which is a typical assupmtion in the scalar literature on this topic, and 2) we allow for serial dependence in the change point model errors, which we assume evolve as a stationary linear process.

This note is organized as follows: In Section \ref{main}, we define the change point model that we consider, and detail the functional CUSUM process as well as the binary segmentation procedure. We also state our main result regarding the consistency of binary segmentation in this setting. Section \ref{sec-proofs} is devoted to the proof of the main result.

\section{Problem Statement and Main Results}\label{main}

Let $\mathcal{H}$ be a separable Hilbert space with inner product denoted $\langle x ,y \rangle$, for $x,y\in \mathcal{H}$, and corresponding norm $\|x\|^2 = \langle x , x \rangle$. Generally we are interested in the case when $\mathcal{H}$ is infinite dimensional, although the below results will still hold when $\mathcal{H}$ is finite dimensional. The example when $\mathcal{H}= L^2([0,1],\mathbb{R})$, where $L^2([0,1],\mathbb{R})$ denotes the space of square integrable real valued functions on the unit interval, is a typical setting for functional data analysis where the data can be thought of as realizations of stochastic processes indexed by some compact set that have square integrable sample paths. Suppose that we observe data $X_1$,$X_2$,...,$X_n$ taking values in $\mathcal{H}$ sequentially, perhaps as a functional time series. All random functions considered in this paper are assumed to be defined over a probability space $(\Omega, \mathcal{F}, Pr)$ and are $\mathcal{F} - \mathcal{B}_{\mathcal{H}}$ measurable, where $\mathcal{B}_{\mathcal{H}}$ denotes the Borel $\sigma$-algebra of subsets of $\mathcal{H}$. In the following discussion we assume $E\|X_i \|<\infty$, so that the mean element of $X_i$ is well defined, see Section 1.5 of \cite{bosq:2000}. We consider in this note the location-error model for the observations $X_i$,
\begin{equation}\label{model}
    X_i = u_i + \varepsilon_i \text{, } \hspace{1cm} i \in \mathbb{Z},
\end{equation}
where $\mathbb{E}[\varepsilon_i] = 0$ for all $i$, and $u_i \in \mathcal{H}$ is the deterministic mean function of $X_i$. We assume that the mean functions follow a generic change point model:

\begin{equation}\label{cp-model}
\mathbb{E}[X_i] = u_i := \mu_j, \hspace{0.1cm} \text{if } v_j^{(n)} < i \leq v_{j+1}^{(n)}, ~~ j \in \{0, 1, \ldots, m_n\}.
\end{equation}

The model \eqref{cp-model} basically specifies that the mean of the series $X_1,...,X_n$ changes at $m_n$ points, $v_1^{(n)},...,v_{m_n}^{(n)}$, with $\mu_1,...,\mu_{m_n}$ denoting the constant means between each change point. Here we allow for the case that $m_n$ increases with the sample size $n$, and also locations of the change points $v_1^{(n)},...,v_{m_n}^{(n)}$ may evolve with $n$. Below we suppress the dependence on $n$ of the change points, and denote them $v_1,...,v_{m_n}$.

The basic goal then is to infer based on the data the number of change points $m_n$ and their locations $v_1,...,v_{m_n}$. A standard method to do this is to employ binary segmentation, which involves sequentially splitting the original sample into two sub-samples based on an initial change point estimate, estimating change points on each sub-sample, and then repeating until some stopping criterion is satisfied. To formulate this method more rigourously, suppose that we have arrived at some point in the procedure at a sub-sample with starting index $l$ and ending index $u$ satisfying $0 < l < u \le n$. In order to identify and estimate change points, one typically considers sequential estimates of the mean function based on the partial sum process $S_k = \sum_{j=1}^{k} X_j$, $1 \leq k \leq n$. Then, in order to estimate changes points on the sub-sample with indices between $l$ and $u$, we consider the generalized CUSUM process $\mathcal{S}_{l,u}^{k}$ defined as
\begin{equation}\label{cumsum}
    \mathcal{S}_{l, u}^{k} = \sqrt{\frac{u-l}{(u-k)(k-l)}}\left[S_k - S_l - \frac{k-l}{u-l}(S_u - S_l)\right].
\end{equation}

Intuitively if there exists one or more change points in the sub-sample with indices between $l$ and $u$, $\|\mathcal{S}_{l, u}^{k}\|$ will be large, and further the point $\hat{k}_{l,u}= \arg\max_{l \le k \le u} \| \mathcal{S}_{l, u}^{k}\|$ estimates a change point. Deciding whether to include or exclude $\hat{k}_{l,u}$ as a potential change point can be determined by the magnitude of $\|\mathcal{S}_{l, u}^{k}\|$: if this magnitude is large enough, say exceeding some threshold $\xi_n$, we then include $\hat{k}_{l,u}$ among the estimated change points, and further segment the sub-sample based on this estimate. The binary segmentation procedure may then be formally described with the following algorithm, as presented in \cite{fryzl:2014:wildbin}.

\floatname{algorithm}{Function}
\renewcommand{\thealgorithm}{}
\begin{algorithm}\label{alg-1}
\caption{BINSEG($l, u, \xi_n$)}
\begin{algorithmic}
\STATE $l \leftarrow$ start index
\STATE $u \leftarrow$ end index
\STATE $\xi_n \leftarrow$ threshold
\IF {$u-l \leq 1$}
        \STATE STOP
\ELSE
        \STATE $k_0 := \underset{l < k < u}{argmax} ||\mathcal{S}_{l, u}^{k}||$
        \STATE $\mathcal{S} = ||\mathcal{S}_{l, u}^{k_0}||$
        \IF {$\mathcal{S} > \xi_n$}
                \STATE add $k_0$ to the set of estimated change points.
                \STATE BINSEG(l, $k_0$, $\xi_n$)
                \STATE BINSEG($k_0$, u, $\xi_n$)
        \ELSE
                \STATE STOP

        \ENDIF
\ENDIF
\STATE
\end{algorithmic}
\end{algorithm}

BINSEG(1, n, $\xi_n$) then returns a set of estimated change points $\hat{\mathcal{V}} = \{\hat{v}_1,...,\hat{v}_{\hat{m}_n}\}$, and an estimated number of change points $\hat{m}_n = |\hat{\mathcal{V}}|.$ We wish to establish the consistency of these estimates under general conditions on the model \eqref{model} and \eqref{cp-model}. As such, we make the following assumptions:

\begin{assumption}\label{err_assump}
Let $\{ w_j, \;\; j \in \mathbb{Z}\}$ be a strong white noise in $\mathcal{H}$ (see Definition 3.1 in \cite{bosq:2000}). The errors $\varepsilon_i$ follow an invertible linear process in $\mathcal{H}$ such that
\begin{equation}\label{err_p}
    \varepsilon_i = \sum_{j=0}^{\infty}h_j(w_{i-j}),
\end{equation}
where $h_0$ is the identity operator, and  $h_j$ is a sequence of bounded linear operators on $\mathcal{H}$ satisfying  $||h_j||_{\mathcal{L}} < a \rho^j$ for some scalars $a>0$ and $\rho \in (0,1)$, with $||\cdot||_{\mathcal{L}}$ denoting the operator norm. Moreover, there exists a $\iota > 0$ such that $\mathbb{E}[\exp (\iota||\varepsilon_0||)] < \infty$.
\end{assumption}

One way in which the results here differ from previous results on the consistency of binary segmentation in this direction in the scalar setting is that we do not assume that the errors are independent and Gaussian, but rather we assume simply that the norms of the errors have sub-exponential tails. Assumption \ref{err_assump} also allows the errors to be serially dependent, and covers many basic functional time series models, including functional autoregressive models under standard conditions.

\begin{assumption}\label{jump_assump}
The minimal jump magnitudes $\Delta = \underset{1 \leq j \leq m_n}{\text{inf}}||\mu_{j + 1} - \mu_{j}||> \eta > 0$ for all $n$.
\end{assumption}

\begin{assumption}\label{dist_assump}
$\underset{0\leq i \leq m_n}{inf}\Bigl(v_{i+1} - v_{i}\Bigr) \geq \delta_n = c_1 n^{1 - \omega}$ where $ 0 \leq \omega < \frac{1}{8}$ and $c_1$ is a positive scalar. \end{assumption}
Assumption \ref{dist_assump} simultaneously restricts the minimal distance between change points and the maximal possible number of changes points, since under this assumption $m_n$ cannot approach infinity faster than $n/\delta_n$, which is on the order of $n^\omega$. This rate is similar to that considered in \cite{venkatraman:1992}. The case when $\omega=0$ corresponds to a finite number of change points. 

\begin{assumption}\label{mean_assump}
$\underset{1 \leq i\leq m+1}{\max} ||\mu_i|| < B < \infty$ for all $n$ and some positive constant $B$.
\end{assumption}

\begin{theorem}\label{main-thm}
Suppose Assumptions \ref{err_assump}--\ref{mean_assump} hold, and let $\hat{m}_n$ and $\hat{v}_1, \hat{v}_2, \ldots, \hat{v}_{\hat{m}_n}$ be respectively the estimated number and locations of the change points, with the estimated change points sorted into increasing order. If the threshold $\xi_n$ in the binary segmentation algorithm satisfies that  $n^{3/8 + e}<\xi_n < n^{1/2 - \omega - e} $ for any positive constant $e< 1/8-\omega$, then  $Pr(\mathcal{B}_n) \to 1$ as $n\to \infty$, where $$\mathcal{B}_n = \{\hat{m}_n = m_n, \;\; |\hat{v}_j - v_j| \le f_n \},$$

and $f_n=n^{5/8 + \omega}\log(n)$. 
\end{theorem}

Theorem \ref{main-thm} establishes that, so long as the threshold $\xi_n$ is chosen to increase with the sample size at an appropriate rate depending on the number and spacing of change points, then  with probability approaching one, the binary segmentation procedure is able to both consistently estimate the number of change points, and localize them to a neighborhood that is both relatively negligible compared to the sample size, as well as the minimal distance between change points, $\delta_n$, defined in Assumption \ref{dist_assump}. When the number of change points is fixed, i.e. $\omega=0$, then consistency can still be achieved when the threshold increases only logarithmically with the sample size, but we do not consider that result here for the sake of brevity.

\section{Proof of Theorem \ref{main-thm}}\label{sec-proofs}

Below we suppose that Assumptions \ref{err_assump}--\ref{mean_assump} are always satisfied. We use $C_i$ to denote unimportant absolute numerical constants. Let $M_k = \sum_{j=1}^{k} u_j$ and $E_k = \sum_{j=1}^{k} \varepsilon_j$. Then we have $\mathcal{S}_{l, u}^{k} = \Theta_{l, u}^{k} + \mathcal{W}_{l, u}^{k}$, where
\begin{equation}\label{signal}
    \Theta_{l, u}^{k} = \sqrt{\frac{u-l}{(u-k)(k-l)}}[M_k - M_l - \frac{k-l}{u-l}(M_u - M_l)]
\end{equation}
and
\begin{equation}\label{error}
    \mathcal{W}_{l, u}^{k} = \sqrt{\frac{u-l}{(u-k)(k-l)}}[E_k - E_l - \frac{k-l}{u-l}(E_u - E_l)]
\end{equation}

The binary segmentation algorithm involves estimating change points on sub-samples, and so in the below Lemmas we generally use $l$ and $u$ to denote the starting and ending index of the sub-sample under consideration.  For any starting index $l$ and ending index $u$, if there are change points between $l$ and $u$, we use the notation $i_0$ and $\beta \geq 0$ to describe the starting index and number of change points between $l$ and $u$, so that
\begin{equation}\label{s_1}
    v_{i_0} \leq l < v_{i_0 + 1} < v_{i_0 + 2} < \ldots < v_{i_0 + \beta} < u \leq v_{i_0 + \beta + 1},
\end{equation}
where $i_0 \leq m_n- \beta$. Let $\mathcal{I} = \{1, 2, \ldots, \beta\}$ be the index set of the change points between $l$ and $u$. The mean function between $v_{i_0 + \gamma}$ and $v_{i_0 + \gamma + 1}$ will be $\mu_{i_0 + \gamma + 1}$ unless otherwise specified for $\gamma \in \mathcal{I}$.\\

\begin{lemma}\label{Err_lemma}
Let $p>1$. Then $Pr(\mathcal{A}_n) \to 1$ as $n\to \infty$, where
\begin{equation}\label{A_n}
    \mathcal{A}_n = \left\{  \underset{1 \leq l < k < u \leq n}{\max}||\mathcal{W}_{l, u}^{k}|| \leq  C_1\log^{p}(n),\;\; \mbox{for all $1 \le l < u \le n$} \right\},
\end{equation}
for some $C_1>0$. 
\end{lemma}
\begin{proof}
Notice that $\mathcal{W}_{l, u}^{k} = \sqrt{\frac{u-l}{(u-k)(k-l)}} \left (\frac{u-k}{u-l}\sum_{j=l+1}^{k} \varepsilon_j - \frac{k-l}{u-l}\sum_{j=k+1}^{u} \varepsilon_{j} \right )$. Therefore,
\begin{IEEEeqnarray}{rCl}
||\mathcal{W}_{l,u}^{k}|| &\leq& \sqrt{\frac{u-k}{u-l}}\biggl | \biggl| \frac{1}{\sqrt{k-l}}\sum_{j=l+1}^{k} \varepsilon_j \biggr| \biggr| + \sqrt{\frac{k-l}{u-l}}\biggl | \biggl |\frac{1}{\sqrt{u-k}}\sum_{j=k+1}^{u} \varepsilon_j\biggr | \biggr |\nonumber\\
&\leq& \biggl | \biggl| \frac{1}{\sqrt{k-l}}\sum_{j=l+1}^{k} \varepsilon_j \biggr| \biggr| + \biggl | \biggl |\frac{1}{\sqrt{u-k}}\sum_{j=k+1}^{u} \varepsilon_j\biggr | \biggr |\nonumber\\
&:=& ||\mathcal{W}_{l, u}^{k, 1}|| + ||\mathcal{W}_{l, u}^{k, 2}||
\end{IEEEeqnarray}
By Bonferroni's inequality,
\begin{IEEEeqnarray}{rCl}
1-Pr(\mathcal{A}_n)&=Pr\left(\underset{l, k, u}{\max} ||\mathcal{W}_{l,u}^{k}|| > C_1\log^{p}(n)\;\; \mbox{for all $1 \le l < u \le n$}\right) \\
&= Pr\left(\underset{l, k, u}{\bigcup}\Bigl|\Bigl|\mathcal{W}_{l,u}^{k}\Bigl|\Bigl| > C_1\log^{p}(n)\right)\nonumber\\
&\leq n^{3}Pr\left(||\mathcal{W}_{l,u}^{k}|| > C_1\log^{p}(n)\right)\nonumber\\
&\leq 2n^{3} \underset{i = 1, 2}{\max}Pr\left(||\mathcal{W}_{l,u}^{k, i}|| > \frac{C_1\log^{p}(n)}{2}\right).\label{bonf_ineq}
\end{IEEEeqnarray}
By Theorem 7.5 of \cite{bosq:2000}, there exist positive constants $C_2$ and $C_3$ so that
\begin{IEEEeqnarray}{rCl}
Pr\left(\Bigl|\Bigl|\mathcal{W}_{l,u}^{k, 1}\Bigr|\Bigr| >\frac{C_1\log^{p}(n)}{2}\right) &=& Pr\left[\frac{1}{\sqrt{k-l}}\Bigl|\Bigl|\mathcal{W}_{l,u}^{k}\Bigl|\Bigl| > \frac{C_1\log^p(n)}{2\sqrt{k-l}}\right]\nonumber\\
&\leq& 4\exp\left(-\frac{(k-l)\frac{C_1^{2}\log^{2p}(n)}{4(k-l)}}{C_2 + C_3 \frac{C_1\log^{p}(n)}{2\sqrt{k-l}}}\right).\label{exp_tail}
\end{IEEEeqnarray}
The same arguments establish this bound for $Pr\left(\Bigl|\Bigl|\mathcal{W}_{l,u}^{k, 2}\Bigr|\Bigr| >{C_1\log^{p}(n)}/{2}\right)$ as well. By substituting (\ref{exp_tail}) back to (\ref{bonf_ineq}), we get that
\begin{IEEEeqnarray}{rCl}
    1-Pr(\mathcal{A}_n) &\leq& 8 \exp\left(3\log(n) - \frac{C_1^{2}\log^{2p}(n)}{4C_2 + 2C_3 \frac{C_1\log^{p}(n)}{\sqrt{k-l}}}\right)\nonumber\\
    &\leq& 8 \exp\left(3\log(n) - \frac{C_1^{2}\log^{2p}(n)}{4C_2 + 2C_1 C_3 \log^{p}(n)}\right).
\end{IEEEeqnarray}
Noticing that 
\begin{equation}
    \frac{C_1^2\log^{2p}(n)}{4C_2 + 2C_1C_3 \log^{p}(n)} \sim \frac{C_1\log^{p}(n)}{2C_3},
\end{equation}
where $a_n \sim b_n$ denotes that $a_n/b_n$ converges to a positive constant, the lemma immediately follows.
\end{proof}
Since the function $\Theta_{l, u}^{k}$ as defined by \eqref{signal} is invariant under translation, we may assume without lose of generality that
\begin{equation}\label{invar}
    M_u - M_l = 0.
\end{equation}
In the following lemmas, when we refer to the mean functions of functional observations $X_k$ within the segment $\{l+1, l+2, \ldots, u\}$, we hence assume the observations have been shifted so that \eqref{invar} holds.
\begin{lemma}\label{qua_lemma}
If $k^{*} = \underset{l < k < u}{\text{argmax}}~||\Theta_{l, u}^{k}||$, then $k^{*} = v_{i_0 + r}$ for some $\gamma \in \mathcal{I}$.
\end{lemma}

\begin{proof}
We consider two cases separately: (i) There is only one change point between $l$ and $u$ and, (ii) There is more than one change point between $l$ and $u$.
Case 1: $|\mathcal{I}| = 1$\\
        There is only one change point between $l$ and $u$, $v_{i_0 + 1}$, which we denote as $v$. Let the mean functions within $[l, v]$ and $[v+1, u]$ be $\mu$ and $\mu '$, respectively. Then we have for $l < k < v$,
        \begin{IEEEeqnarray}{rCl}
            ||\Theta_{l, u}^{k}|| &=&\sqrt{\frac{u-l}{(u-k)(k-l)}} (k-l) ||\mu|| \nonumber\\
            &=& \sqrt{u-l} \sqrt{\frac{k-l}{u-k}} ||\mu|| \label{sp_1}
        \end{IEEEeqnarray}
        From this it is clear that $||\Theta_{l, u}^{k}||$ is either a monotonically increasing or identically zero as a function of $k$ satisfying  $l < k < v$.
        Similarly, for any $k$ such that $v <k < u$, we have
        \begin{IEEEeqnarray}{rCl}
            ||\Theta_{l, u}^{k}|| &=& \sqrt{\frac{u-l}{(u-k)(k-l)}} ||M_k - M_l||\nonumber\\
            &=& \sqrt{\frac{u-l}{(u-k)(k-l)}} ||M_u - M_k||\nonumber\\
            &=& \sqrt{u-l} \sqrt{\frac{u-k}{k-l}}||\mu'|| \label{sp_2}
        \end{IEEEeqnarray}
        which is also monotonically decreasing or identically zero as a function with respect to $k$.
        Moreover, $||\Theta_{l, u}^{k}||$ cannot be zero over the entire segment $[l,u]$ as $||\mu - \mu'|| > \Delta$ implies $\max(||\mu||, ||\mu'||) > \frac{\Delta}{2} > 0$. Thus, $||\Theta_{l, u}^{v}|| = \underset{l < k < u}{\max}||\Theta_{l_u}^{k}||$, and the maximum $v$ is unique. \\

Case 2: $|\mathcal{I}| > 1$\\
        Let $v = v_{i_0 + \alpha}$ and $v' = v_{i_0 + \alpha + 1}$ for some $\alpha \in \mathcal{I}$. If $v$ is the right-most change point, let $v' = u$. The case when $v$ is the left most change point can be handled similarly. Otherwise, let $D_{v, k} = ||\Theta_{l, u}^{v}|| - ||\Theta_{l, u}^{k}||$ and $D_{v', k} = ||\Theta_{l, u}^{v'}|| - ||\Theta_{l, u}^{k}||$.  Suppose $d_j = (v_{i_0 + j} - l)/(u-l)$ for $j \in \mathcal{I}$, and let $d^*_1=(v - l)/(u-l)$, $d^*_2=(v' - l)/(u-l)$. Then $0 =d_0 < d_1 < d_2 < \ldots < d_\beta < d_{\beta+1} = 1$. Therefore, for any $k$ between $v$ and $v'$, we may rewrite $\Theta_{l, u}^{k}$ as follows:
        
        \begin{IEEEeqnarray}{rCl}
        \Theta_{l, u}^{k} &=& \sqrt{\frac{u-l}{(u-k)(k-l)}}(M_k - M_l)\nonumber\\
        &=&\sqrt{u-l}\frac{\sum_{j=1}^{r}(d_j - d_{j-1})\mu_{i_0 + j} + (\frac{k-l}{u-l} - d_r)\mu_{i_0 + r +1}}{\sqrt{\frac{k-l}{u-l}\frac{u-k}{u-l}}}.\nonumber
        \end{IEEEeqnarray}
        Let $x=x(k) =(k-l)/(u-l)$. Then
        \begin{IEEEeqnarray}{rCl}
            ||\Theta_{l, u}^{k}||^2 &=& (u-l)\frac{||\sum_{j=1}^{\alpha}(d_j - d_{j-1})\mu_{i_0 + j} + (x - d_\alpha)\mu_{i_0 + \alpha +1}||^2}{x(1-x)} \nonumber\\
            &:=& (u-l)f(x)
        \end{IEEEeqnarray}
        Let $s(x) = ||\sum_{j=1}^{\alpha}(d_j - d_{j-1})\mu_{i_0 + j} + (x - d_\alpha)\mu_{i_0 + \alpha +1}||^2$, the numerator of $f(x)$. If we further simplify $s$:
        \begin{align}\label{rCl}
            s(x) &= \sum_{j=1}^{\alpha}(d_j - d_{j-1})^2||\mu_{i_0+j}||^2 + (x-d_\alpha)^2||\mu_{i_0 + \alpha + 1}||^2 \nonumber\\
            &+ 2(x - d_\alpha)\sum_{j=1}^{\alpha}(d_j - d_{j-1})\langle\mu_{i_0 + j}, \mu_{i_0 + \alpha + 1}\rangle\nonumber\\
            &=\left(\sum_{j=1}^{\alpha}(d_j - d_{j-1})||\mu_{i_0+j}|| + (x-d_\alpha)||\mu_{i_0 + \alpha + 1}||\right)^2\nonumber\\
            &+ 2(x - d_\alpha)\sum_{j=1}^{\alpha}(d_j - d_{j-1})\Bigl(\langle\mu_{i_0 + j}, \mu_{i_0 + \alpha + 1}\rangle - ||\mu_j||\cdot||\mu_{i_0 + \alpha + 1}||\Bigr)\nonumber\\
            &=\bigl(a'x+b'\bigr)^2 + 2(x-d_\alpha)t,
        \end{align}
        where $a' = ||\mu_{i_0 + \alpha + 1}||$, $b' =\sum_{j=1}^{\alpha}(d_j - d_{j-1})||\mu_{i_0+j}|| - d_\alpha a'$ and \\ 
        $t = \sum_{j=1}^{\alpha}(d_j - d_{j-1})\Bigl(\langle\mu_{i_0 + j}, \mu_{i_0 + \alpha + 1}\rangle - ||\mu_j||\cdot||\mu_{i_0 + \alpha + 1}||\Bigr)$. Notice that by the Cauchy-Schwarz inequality, $t \le 0$. Moreover, \eqref{rCl} can be represented as
        \begin{IEEEeqnarray}{rCl}
            s(x) &=& a'^2 x^2 + 2(t + a'b')x + b'^2-2td_\alpha\nonumber\\
            &:=& ax^2 + bx + c.\nonumber
        \end{IEEEeqnarray}
        If $f(x)$ is extended to the open unit interval as $s(x)/[x(1-x)]^2$, then what we now wish to show is that $f$ achieves a maximum over the interval $[d_1^*,d_2^*]$ at a point on the boundary, at either $d_1^*$ or $d_2^*$. The derivative of $f(x)$ is
        \begin{IEEEeqnarray}{rCl}
            f'(x) &=& \frac{(a+b)x^2 + 2cx - c}{[x(1-x)]^2}=: \frac{g(x)}{[x(1-x)]^2},
        \end{IEEEeqnarray}
        where $g(x)$ is a quadratic function with vertex $-\frac{c}{a+b}$, when $(a+b)\ne0$. First notice that $g(0) = -c = -(b'^2 - 2td_\alpha) \leq 0$. There are three scenarios that we consider: (i) $a + b = 0$. (ii) $a + b > 0$, and (iii) $a + b < 0$.\\
        Scenario 1: $a + b = 0$\\
            First we claim that in this case $c \neq 0$. If $c=0$, then we will have $b' = t = 0$. Notice that,
            \begin{IEEEeqnarray}{rCl}
                a + b &=& a'^2 + 2(t+ a'b')\\ \notag
                &=& a'^2 
            \end{IEEEeqnarray}
            Therefore, $a+b = 0$ and $c= 0$ implies $a' = ||\mu_{i_0 + \alpha + 1}|| = 0$ and $b' = \sum_{j=1}^{\alpha}(d_j - d_{j-1})||\mu_{i_0 + j}|| - 0 = 0$. Combining these lead to $\mu_{i_0 + 1} = \mu_{i_0 + 2} =\ldots = \mu_{i_0 + \alpha + 1} = 0$ which contradicts Assumption \ref{jump_assump}. Under this scenario, $g(x) = c(2x-1)$ and $g(0) = -c < 0$. This implies $f(x)$ decreases on $[0, 0.5]$ and increases on $[0.5, 1]$. \\
            Scenario 2: $a+b > 0$\\
            The vertex of $g(x)$ is $- \frac{c}{a+b}$ which is negative in this case. Therefore, $g(x)$ is negative from $0$ to some real number $x_0$, and is positive from $x_0$ to infinity if $c > 0$. If $c=0$, then $g(x)$ is always positive, which implies $f(x)$ strictly increases on $[0, 1]$.
            Scenario 3: $a+b < 0$\\
            In this scenario, The vertex $-\frac{c}{a+b}$ will be positive and the maximum of $g(x)$ is $-[\frac{c^2}{a+b} + c]$. If the maximum of $g(x)$ is negative, then $g(x)$ is always negative and $f(x)$ will be strictly decreasing. Otherwise, we will have $-[\frac{c^2}{a+b} + c] > 0$ which implies $\frac{c}{a+b} < -1$. The roots of $g(x)$ are $x_1 = - \sqrt{\frac{c}{a+b}(\frac{c}{a+b} + 1)} - \frac{c}{a+b}$ and $x_2 = + \sqrt{\frac{c}{a+b}(\frac{c}{a+b} + 1)} - \frac{c}{a+b}$. Clearly, $x_2 > 1$. Therefore, $g(x)$ is either positive from $0$ to $1$ or negative from $0$ to $x_1$ and positive from $x_1$ to $1$. Once again, $f(x)$ over $[0,1]$ is either strictly increasing, or decreasing and then increasing, respectively.\\

It follows then in all cases that $f(x)$ is maximized over $[d_1^*,d_2^*]$ at either $d_1^*$ or $d_2^*$, from which the lemma follows. 

\end{proof}

\begin{remark}
From equation \eqref{sp_1} and \eqref{sp_2}, $||\Theta_{l, u}^{l}||= ||\Theta_{l, u}^{u}|| = 0$
\end{remark}
The following conditions are imposed through the below Lemmas \ref{max_order} - \ref{dist_lemma} on the sub-sample between $l$ and $u$:
\begin{equation}\label{s_2}
    l < v_{i_0 + r} - \xi \delta_n < v_{i_0 + r} + \xi\delta_n < u ~~ \text{ for some } r\in \mathcal{I} \text{ and } \xi\in \left[\frac{1}{2}, 1\right)
\end{equation}
\begin{equation}\label{s_3}
    \max(\min(v_{i_0 + 1} - l, l - v_{i_0}), \min(u-v_{i_0 + \beta}, v_{i_0 + \beta + 1} - u)) \leq n^{\frac{5}{8} + \omega}\log(n)
\end{equation}

\begin{lemma}\label{max_order}
 $\Theta_{l, u} \geq \frac{\xi\Delta}{2} \frac{\delta_n}{\sqrt{n}} = O(n^{1/2 - \omega})$
\end{lemma}
\begin{proof}
Let $a = \underset{l < k < u}{\max}||M_k - M_l||$. Then we aim to show that $a \geq  \frac{\xi\Delta}{4}\delta_n$ under condition \eqref{s_2} and Assumptions \ref{err_assump}--\ref{dist_assump}. Let $v =v_{i_0 + r}$ and $v' = v_{i_0 + r + 1}$(if $v$ is the right-most change point, let $v' = u$). Further assume $E[X_v] = \mu$ and $E[X_v'] = \mu'$. Since $||\mu - \mu'|| \geq \Delta$, we get by the reverse triangle inequality that $\max\left(||\mu||,||\mu'||\right) \geq \frac{\Delta}{2}$. Moreover, immediately we get from \eqref{s_2} and \eqref{s_3} that there is no change point between $[v - \xi\delta_n, v)$ and $(v, v +\xi\delta_n]$. it leads to
\begin{equation}\label{con-eq}
\max(||M_v - M_{v-\xi\delta_n}||, ||M_{v+ \xi\delta_n} -M_v||) \geq \frac{\xi\Delta}{2}\delta_n.
\end{equation}
Then we claim that $\underset{l < k < u}{\max}||M_k - M_l|| \geq \frac{\xi\Delta}{4}\delta_n$. If not, 
\begin{equation}
    \begin{cases}
        ||M_{v + \xi\delta_n} - M_l|| &< \frac{\xi\Delta}{4}\delta_n, \\
        ||M_{v} - M_l|| &< \frac{\xi\Delta}{4}\delta_n, \;\; \mbox{and}\\
        ||M_{v -  \xi\delta_n} - M_l|| &< \frac{\xi\Delta}{4}\delta_n. \\
    \end{cases}
\end{equation}
These with the triangle inequality contradict \eqref{con-eq}. Furthermore, because $\frac{u-l}{(u-k)(k-l)} \leq \frac{n}{4}$, we have
\begin{IEEEeqnarray}{rCl}
    \Theta_{l, u} &=& \underset{l<k<u}{\max}\frac{1}{\sqrt{\frac{(u-k)(k-l)}{u-l}}}||M_v - M_l||\nonumber\\
    &\geq& \frac{1}{\sqrt{\frac{n}{4}}} \frac{\xi\Delta}{4}\delta_n = \frac{\xi\Delta}{2}\frac{\delta_n}{\sqrt{n}} \nonumber
\end{IEEEeqnarray}
\end{proof}

\begin{lemma}\label{bound_lemma}
Suppose the sub-sample with starting and ending index $l$ and $u$ respectively satisfy \eqref{s_2} and \eqref{s_3}. Under A\ref{err_assump} - A\ref{mean_assump}, for those $C_1$ and $p$ in Lemma \ref{Err_lemma}, let $\Theta_{l, u} = \underset{l<k<u}{\max}||\Theta_{l, u}^{k}||$ and $v$ be a change point that satisfies $l<v<u$,
\begin{equation}\label{cond}
    \Bigl|\Bigl |\Theta_{l, u}^{v}\Bigr| \Bigr | > \Theta_{l, u} - 2C_1 \log^{p}(n)
\end{equation}
Then we have that
\begin{equation}\label{bound_signal}
    C_l \frac{\delta_n}{\sqrt{n}}\leq \Bigl|\Bigl |\Theta_{l, u}^{v}\Bigr| \Bigr | \leq C_u \frac{n}{\sqrt{\delta_n}}
\end{equation}
for some positive real numbers $C_l$ and $C_u$, independent from sample size.
\end{lemma}
\begin{proof}
Let $v=v_{i_0 + i}$ for some $i \in \mathcal{I}$.
First notice that either $\min(v-l, u-v) \leq n^{5/8 + \omega}\log(n)$ or $\min(v-l, u-v) \geq \delta_n - n^{5/8+ \omega}\log(n)$ by the setting \eqref{s_1} and \eqref{s_2}. If the former holds, suppose $v-l \leq n^{5/8 + \omega}\log(n)$. Immediately we get $M_v - M_l = (v-l)\mu$ by assuming the mean function between $l$ and $v$ is $\mu$. This leads to $ (v-l)||\mu|| \leq B n^{5/8+ \omega}\log(n)$. However, $n^{5/8+ \omega}\log(n) = o(\delta_n)$ as $\omega < 1/8$, which with \eqref{cond} contradicts Lemma \ref{max_order}. If instead $u-v \leq n^{5/8+ \omega}\log(n)$ we may obtain a similar contradiction. Thus,
\begin{equation}\label{min_d}
    \min(v-l, u-v) \geq \delta_n - n^{5/8+ \omega}\log(n)
\end{equation}
Since $M_v - M_l =(v_{i_0 +1} - l)\mu_{i_0 + 1} +  (v_{i_0 +2} - v_{i_0 +1})\mu_{i_0 + 2} + \ldots + (v_{i_0 + r} - v_{i_0 + r -1}) \mu_{i_0 + r}$, a linear combination of functions in $\mathcal{H}$, we  define $\mu^{*} = (M_v - M_l)((v_{i_0 + r} - l))$, which is an  $\mathcal{H}$, and $||\mu^{*}|| \leq B$. Further, $||\mu^{*}|| \neq 0$ as $||M_v - M_l|| \neq 0$ according to \eqref{cond} and Lemma \ref{max_order}. Therefore,
\begin{IEEEeqnarray}{rCl}
    ||\Theta_{l, u}^{v}|| &=& \sqrt{\frac{u-l}{(v-l)(u-v)}}||M_v-M_l||\nonumber\\
    &=& \sqrt{\frac{u-l}{u-v}}\sqrt{v-l}||\mu^{*}||\nonumber\\
    &\leq& \sqrt{\frac{n}{\delta_n - n^{5/8+ \omega}\log(n)}}\sqrt{n}B\nonumber\\
    &\leq& C_4B \frac{n}{\sqrt{\delta_n}} \label{ineq_1}\\
    &=& C_5 \frac{n}{\sqrt{\delta_n}},\nonumber
\end{IEEEeqnarray}
for some constants $C_4$ and $C_5$. This follows since $n^{5/8+ \omega}\log(n) = o(\delta_n)$, and hence
\begin{equation}
    \underset{n\to \infty}{\lim} \frac{\delta_n - n^{5/8+ \omega}\log(n)}{\delta_n} =1
\end{equation}
Therefore, for any $\epsilon >0$, there exists an integer $K$ such that $\delta_n - n^{5/8+ \omega}\log(n) \le (1+\epsilon) \delta_n$, for all $n\ge K$. So we may take $C_5 = \sqrt{1 + \epsilon}$ and the upper bound follows.

The lower bound for $||\Theta_{l, u}^{v}||$ is obtained by \eqref{cond} and Lemma \ref{max_order}. Indeed,
\begin{IEEEeqnarray}{rCl}
    ||\Theta_{l, u}^{v}|| &>& \Theta_{l, u} - 2\lambda\log^{1/p}(n)\nonumber\\
    &\geq& \frac{\xi\Delta}{2} \frac{\delta_n}{\sqrt{n}}- 2\lambda\log^{1/p}(n)\nonumber\\
    &\geq& C_l \frac{\delta_n}{\sqrt{n}} = O(n^{1/2 - \omega}). \label{ineq_2}
\end{IEEEeqnarray}
Equation \eqref{ineq_2} uses the same argument used to establish the upper bound for $||\Theta_{l, u}^{v}||$.
\end{proof}

\begin{lemma}\label{diff_lemma}
Suppose $l, u$ satisfy \eqref{s_2} and \eqref{s_3}. Let $v$ the be change point satisfying \eqref{cond}. Then there exists $\gamma_1$ and $\gamma_2$ satisfying $0< \gamma_i \leq n^{5/8 + \omega}\log(n)$, $i=1,2$ so that
\begin{equation}
    ||\Theta_{l, u}^{v + \gamma_1}|| < ||\Theta_{l, u}^{v}|| - 2C_1 \log^{p}(n)
\end{equation}
\begin{equation}
    ||\Theta_{l, u}^{v - \gamma_2}|| < ||\Theta_{l, u}^{v}|| - 2C_1 \log^{p}(n)
\end{equation}

\end{lemma}
\begin{proof}
The result follows directly from Lemma 2.6 in \cite{venkatraman:1992}. Indeed,
\begin{equation}
    ||\Theta_{l, u}^{v}|| = \sqrt{\frac{u-l}{(u-v)(v-l)}}||M_v - M_l|| = \sqrt{\frac{u-l}{(u-v)(v-l)}} a
\end{equation}
which share the same form as the quotient in Lemma 2.6 of \cite{venkatraman:1992}.
\end{proof}

\begin{lemma}\label{dist_lemma}
Suppose $l, u$ satisfies condition \eqref{s_2} and \eqref{s_3}. Then the estimated location of the change point $\hat{k} = \underset{l<k<u}{\text{argmax}} ~||\mathcal{S}_{l, u}^{k}||$ satisfies $$Pr\left(\biggl\{ dist(\hat{k}, \mathcal{V}) > f_n \biggr\} \cap \mathcal{A}_n\right) \to 0, ~\text{ as } n \to \infty,$$ where $\mathcal{V} = \{v_{i_0 + 1}, v_{i_0 + 2}, \ldots, v_{i_0 + \beta}\}$, $f_n = n^{5/8 + \omega}\log(n) \text{, and } \text{dist}(x, A) = inf\{|x- a|~|~ \forall a \in A\}.$
\end{lemma}

\begin{proof}
On $\mathcal{A}_n$, we have that $\underset{l<k<u}{\max}||\mathcal{W}_{l, u}^{k}|| \leq C_1 \log^{p}(n)$. Let $k^* = \underset{l <k<u}{\text{argmax}}||\Theta_{l, u}^{k}||$. Therefore, on $\mathcal{A}_n$, we have by the triangle inequality that
\begin{IEEEeqnarray}{rCl}
    \Theta_{l, u} = ||\Theta_{l, u}^{k^{*}}|| &=& ||\Theta_{l, u}^{k^*} + \mathcal{W}_{l, u}^{k^*} - \mathcal{W}_{l, u}^{k^*}|| \nonumber \\
    &\leq& ||\mathcal{S}_{l, u}^{k^*}|| + C_1 \log^{p}(n)\nonumber\\
    &\leq& ||\mathcal{S}_{l, u}^{\hat{k}}|| + C_1 \log^{p}(n)\nonumber\\
    &\leq& ||\Theta_{l, u}^{\hat{k}}|| + 2C_1 \log^{p}(n)\nonumber
\end{IEEEeqnarray}
We aim now to show that 

$$
\biggl\{ dist(\hat{k}, \mathcal{V}) > f_n \biggr\} \cap \mathcal{A}_n
$$
is empty for all $n$ sufficiently large. Suppose the above set is non-empty and contains an element $\nu$. Regardless of the value of $\hat{k}$, it is between $v := v_{i_0 + r}$ and $v' := v_{i_0 + r + 1}$ for some $r \in \mathcal{I}$. If $r = \beta$, then let $v' = u$. Lemma \ref{qua_lemma} and the above inequalities imply that $||\Theta_{l, u}^{k}||$ is either monotonic or decreases and then increases between $v$ and $v'$, and so
\begin{equation*}
    \max(||\Theta_{l,u}^{v}||, ||\Theta_{l, u}^{v'}||) > ||\Theta_{l,u}^{\hat{k}}|| > \Theta_{l, u} - 2C_1 \log^{p}(n).
\end{equation*}
If $v$ is the right-most change point between $l$ and $u$, then $v-u \geq \delta_n - f_n$ by \eqref{min_d} as $||\Theta_{l, u}^{u}|| = 0$. Note that on the sample point $\nu$, $\hat{k} \in [v + f_n, v' - f_n]$.\\
Suppose that $||\Theta_{l,u}^{v}|| > ||\Theta_{l,u}^{v'}||$. Then again, according to the proof of Lemma \ref{qua_lemma},  $||\Theta_{l, u}^{k}||$ is either strictly decreasing or decreasing then increasing from $v$ to $v'$. Since $||\Theta_{l,u}^{v}||$ satisfies \eqref{cond}, there exists a $\gamma \in (0, f_n)$ by Lemma \ref{diff_lemma} so that
\begin{equation*}
    ||\Theta_{l,u}^{v+ \gamma}|| < ||\Theta_{l,u}^{v}|| - 2C_1 \log^{p}(n).
\end{equation*}
 for all $n$ sufficiently large. If $||\Theta_{l, u}^{k}||$ is strictly decreasing from $v$ to $v'$, then
\begin{IEEEeqnarray}{rCl}
    \Theta_{l, u} \geq ||\Theta_{l, u}^{v}|| &>& ||\Theta_{l,u}^{v+ \gamma}|| + 2C_1 \log^{p}(n) \nonumber\\
    &>& ||\Theta_{l,u}^{\hat{k}}||+ 2C_1 \log^{p}(n)\nonumber\\
    &\geq& \Theta_{l, u},
\end{IEEEeqnarray}
which is a contradiction. Therefore, $||\Theta_{l,u}^{k}||$ must be first decreasing and then increasing, and further must be increasing for all integers $k$ such that $k > \hat{k}$. In this case, we will have $||\Theta_{l, u}^{v'}|| > ||\Theta_{l, u}^{\hat{k}}||$ which implies $||\Theta_{l, u}^{v'}||$ satisfies \eqref{cond}. So again by Lemma \ref{diff_lemma},
\begin{equation*}
    ||\Theta_{l, u}^{v' - \gamma'}|| < \Theta_{l, u}^{v'} - 2C_1 \log^{p}(n),
\end{equation*}
for some $r'\in (0, f_n)$ and for all $n$ sufficiently large. The same argument could be applied here
\begin{IEEEeqnarray}{rCl}
    \Theta_{l, u} \geq ||\Theta_{l, u}^{v'}|| &>& ||\Theta_{l,u}^{v' - \gamma}|| + 2C_1 \log^{p}(n)\nonumber\\
    &>& ||\Theta_{l,u}^{\hat{k}}||+ 2C_1 \log^{p}(n)\nonumber\\
    &\geq& \Theta_{l, u},
\end{IEEEeqnarray}
giving another contradiction. If $||\Theta_{l,u}^{v}|| < ||\Theta_{l,u}^{v'}||$, a similar argument can be applied. In summary, the above argument shows the event $$\biggl\{ dist(\hat{k}, \mathcal{V}) > f_n \biggr\} \cap \mathcal{A}_n$$ is empty for all $n$ sufficiently large, and so the result follows by the continuity of the probability measure.
\end{proof}

\begin{lemma}\label{order_no_break}
Suppose $l$ and $u$ such that one of the following conditions are satisfied with $f_n$ defined in Lemma \ref{dist_lemma}
\begin{enumerate}[(i)]
    \item $\beta = 0$, $v_{i_0} < l < u < v_{i_0 + 1}$
    \item  $\beta = 1$, $\min(v_{i_0 + 1} - l, u - v_{i_0 + 1}) \leq f_n $\\
    \item $\beta = 2$, $\max(v_{i_0 + 1} - l, u - v_{i_0 + 2}) \leq f_n $
\end{enumerate}
Then on the set $\mathcal{A}_n$, then for large n, we will have
\begin{equation}
    \underset{l < k < u}{\max}||\mathcal{S}_{l,u}^{k}|| \leq C_6 n^{3/8}\sqrt{\log(n)}
\end{equation}
where $C_6$ is independent from sample size.
\end{lemma}
\begin{proof}
Indeed,
\begin{enumerate}[start=1,label={\bf Case\arabic* :}, leftmargin=7em]
    \item $\beta = 0$\\
        If the first condition holds, then there is no change point between $l$ and $u$. Therefore, $||\mathcal{S}_{l, u}^{k}|| = ||\mathcal{W}_{l, u}^{k}||$. The uniform bound then will be $C_1 \log^{p}(n)$ as $\pmb{\varepsilon} \in \mathcal{A}_n$
    \item $\beta = 1$\\
        If the second condition holds then
        \begin{IEEEeqnarray}{rCl}
            \Theta_{l, u} &=&||\Theta_{l, u}^{v_{i_0 + 1}}|| = \sqrt{u-l}\sqrt{\frac{v_{i_0 + 1} - l}{u-v_{i_0 + 1}}}||\mu||\nonumber\\
            &\leq& \sqrt{2}B \sqrt{\min(v_{i_0 + 1}-l, u-v_{i_0 + 1})}\nonumber\\
            &\leq& \sqrt{2}B \sqrt{n^{5/8 + \omega}\log(n)}\nonumber\\
            &\leq& C_1 n^{3/8}\sqrt{\log(n)} \nonumber
        \end{IEEEeqnarray}
        Therefore,
        \begin{IEEEeqnarray}{rCl}
            \underset{l<k<u}{\max}||\mathcal{S}_{l, u}|| &\leq& \underset{l<k<u}{\max}||\Theta_{l, u}|| + \underset{l<k<u}{\max}||\mathcal{W}_{l, u}||\nonumber\\
            &\leq& \Theta_{l, u} + C_1 \log^{p}(n)\nonumber\\
            &\leq& C_1^{*} n^{3/8}\sqrt{\log(n)} ~~~~\text{ as } C_1 \log^{p}(n) = o\Bigl( n^{3/8}\sqrt{\log(n)}\Bigr) \nonumber
        \end{IEEEeqnarray}
    \item $\beta = 2$\\
        If the third condition holds then
        \begin{IEEEeqnarray}{rCl}
            \Theta_{l, u} &=& \max(||\Theta_{l, u}^{v_{i_0 + 1}}||, ||\Theta_{l, u}^{v_{i_0 + 2}}||)\nonumber\\
            &\leq& \sqrt{2}B \sqrt{\max(v_{i_0 + 1} - l, u-v_{i_0 + 2})}\nonumber\\
            &\leq& \sqrt{2}B \sqrt{n^{5/8 + \omega}\log(n)}\nonumber\\
            &\leq& C_2 n^{3/8}\sqrt{\log(n)} \nonumber
        \end{IEEEeqnarray}
        Apply the same argument as the second case, we will have
        \begin{equation}
            \underset{l<k<u}{\max}||\mathcal{S}_{l, u}|| \leq C_2^{*} n^{3/8}\sqrt{\log(n)}
        \end{equation}

\end{enumerate}
Therefore, we will have
\begin{equation}
    \Theta_{l, u} \leq \max(C_1^{*}, C_2^{*}) n^{3/8}\sqrt{\log(n)} = C_6 n^{3/8}\sqrt{\log(n)}
\end{equation}
\end{proof}

\begin{proof}[Proof of Theorem \ref{main-thm}] First notice that $Pr[\mathcal{B}_n^{c}] = Pr[\mathcal{B}_n^{c}\cap \mathcal{A}_n] + Pr[\mathcal{B}_n^{c}\cap \mathcal{A}_n^{c}]$. By Lemma \ref{Err_lemma}, $ Pr[\mathcal{B}_n^{c}\cap \mathcal{A}_n^{c}] \to 0$ as $n\to \infty$. On $\mathcal{A}_n$, the binary segmentation procedure begins by letting $l=0$ and $u=n$. $l$ and $u$
satisfy conditions \eqref{s_1}, \eqref{s_2} and \eqref{s_3} as long as $m_n \ge 1$. Thus,  for all $n$ sufficiently large and on the set $\mathcal{A}_n$ the first
estimator $\hat{k}_1$ will fall into a neighbourhood of width $f_n$ of some change point $v_{i_0 + r}$ for some $r \in \{1, 2, \ldots, m_n\}$ by Lemma \ref{dist_lemma}. Binary segmentation starts again on the new sub segments with starting and ending indices $l_1^{1}=0$, $u_1^{1} = \hat{k}_1$, and start $l_2^{1} = \hat{k}_1$, $u_2^{1} = n$. The procedure is repeated on each segment as long as \eqref{s_1} and \eqref{s_3} are satisfied on each segment. If we have detected less than $m_n$ change points, then there must exist a segment $\{l, l+1, \ldots, u\}$ such that \eqref{s_1} and \eqref{s_2} hold. Moreover, $l$ and $u$ are estimated location of change point which implies, on $\mathcal{A}_n$ for $n$ sufficiently large, $l$ and $u$ must be within $f_n$ of one of the true change points, which implies \eqref{s_3} is satisfied. Thus, one more change point in that segments will be detected by Lemma \ref{max_order}. Hence, $\hat{m}_n \geq m_n$ necessarily on $\mathcal{A}_n$ for $n$ sufficiently large. Once we have detected $m_n$ change points, the end-points of each segment will satisfy one of the cases in Lemma \ref{order_no_break}. It follows then that $\mathcal{B}_n^{c}\cap \mathcal{A}_n$ is the empty for all $n$ sufficiently large, and hence $Pr[\mathcal{B}_n^{c}\cap \mathcal{A}_n]\to 0$ as $n\to \infty$ by the continuity of the probability measure. Thus, with the above, $Pr[\mathcal{B}_n^{c}] \to 0$ as $n \to \infty$.
\end{proof}

\bibliographystyle{apalike}
\bibliography{Greg_bib}

\end{document}